\documentclass[12pt]{article}
\usepackage{amsmath}
\usepackage{graphicx}
\usepackage{amsfonts}
\usepackage{amssymb}
\usepackage{stmaryrd}
\usepackage{amscd}
\newcommand{\bc}{\begin{center}}
\newcommand{\ec}{\end{center}}
\newcommand{\be}{\begin{equation}}
\newcommand{\ee}{\end{equation}}
\newcommand{\bea}{\begin{eqnarray}}
\newcommand{\eea}{\end{eqnarray}}
\newcommand{\ba}{\begin{array}}
\newcommand{\ea}{\end{array}}
\newcommand{\edc}{\end{document}}

\begin{document}
\thispagestyle{empty}
\begin{center}

{\textbf{$\ell$- VOLTERRA QUADRATIC STOCHASTIC OPERATORS:\\ LYAPUNOV FUNCTIONS, TRAJECTORIES}}\\[4mm]
{\bf U. A. Rozikov}$^{1,2}$ and {\bf A. Zada}$^2$\\[2mm]
$^1$Institute of Mathematics and Information Technologies,\\ Tashkent, Uzbekistan.\\
email: rozikovu@yandex.ru \\
$^2$Abdus Salam School of Mathematical Sciences, GCU,\\
35-C-2, Gulberg III, Lahore, Pakistan.\\
email: zadababo@yahoo.com\\
\vspace{0.3cm}
\end{center}
{\bf Abstract.} We consider $\ell$-Volterra quadratic stochastic operators defined on
$(m-1)$-dimensional simplex, where $\ell\in\{0,1,...,m\}$. Under some conditions on coefficients of
such operators we describe Lyapunov functions and apply them to obtain upper estimates for the set
of $\omega$- limit points of trajectories. We describe a set of fixed points of $\ell$-Volterra
operators. {\vskip .3cm \noindent


\section{Introduction}
A quadratic stochastic operator (QSO) has meaning of a population evolution operator (see [6]-[8]),
which arises as follows. Consider a population consisting of $m$ species. Let $x^{0} =
(x_{1}^{0},...,x_{m}^{0})$ be the probability distribution of species in the initial generations,
and $P_{ij,k}$ the probability that individuals in the $i$th and $j$th species interbreed to
produce an individual $k$. Then the probability distribution $x'= (x_{1}',...,x_{m}')$ (the state)
of the species in the first generation can be found by the total probability i.e.
$$ x'_k=\mathop {\sum} \limits^{m}_{i,j=1}P_{ij,k}x^{0}_{i}x^{0}_{j} ,\,\,\,\,k=
1,...,m. \eqno(1)$$
 This means that the association $x^{0}\rightarrow x'$
defines a map $V$ called the evolution operator. The population
evolves by starting from an arbitrary state $x^{0}$, then passing to
the state $x'= V(x)$ (in the next "generation"), then to the state
$x''=V(V(x))$, and so on.
 Thus states of the population described by the following dynamical
 system
         $$ x^{0},\ \ x'= V(x), \ \ x''=V^{2}(x),\ \  x'''= V^{3}(x),...$$
Note that $V$ (defined by (1)) is a non linear (quadratic) operator, and it is higher dimensional
if $m\geq 3$. Higher dimensional dynamical systems are important but there are relatively few
dynamical phenomena that are currently understood ([1],[2],[9]).

In [11] we considered a class of nonlinear (quadratic) operators which is called $\ell$-Volterra
operators and the difference of $\ell$- Volterra quadratic operators from known quadratic operators
are discussed. Some invariant (in particular some fixed points) sets for $\ell$-Volterra operators
are described. Also we described a family of $\ell$-Volterra operators each element
 of which has cyclic orbits generated by several
 vertices of the simplex. It is shown that the set of all
 $\ell$-Volterra operators is convex, compact and its extremal points are constructed.
 For 1-Volterra operators and 2-Volterra
 operators defined on a two dimensional simplex the limit behavior of all trajectories
(orbits) are studied.

In this paper we continue the investigations of $\ell$-Volterra quadratic operators. Under some
conditions on coefficients of such operators we describe Lyapunov functions and apply them to
obtain upper estimates for the set of $\omega$- limit points of trajectories. We describe a set of
fixed points of the $\ell$-Volterra operators. This paper also contains many remarks with
comparisons of $\ell$-Volterra operators and Volterra ones.

\section{Definitions}

The quadratic stochastic operator (QSO) is a mapping of the simplex.

$$ S^{m-1}=\left\{x=(x_1,...,x_m)\in {\bf R}^m: x_i\geq
0,\,\sum^m_{i=1}x_i=1 \right\} \eqno(2)$$ into itself , of the form
$$
    V: x'_k= \sum^m_{i,j=1}P_{ij,k}x_ix_j,\,\,\,k=
    1,...,m,\eqno(3)$$
where $P_{ij,k}$ are coefficients of heredity and
$$P_{ij,k}\geq 0,\ \ P_{ij,k}= P_{ji,k},\ \ \sum^m_{k=1}P_{ij,k}= 1, (i,j,k=
1,....,m).\eqno(4)$$

 Thus each quadratic stochastic operator $V$ can
be uniquely defined by a cubic matrix
$\textbf{P}=\left(P_{ij,k}\right)^n _{i,j,k=1}$ with conditions
(4).

Note that each element $x\in S^{m-1}$ is a probability distribution
on $E= \{1,...,m\}$. The population evolves by starting from an
arbitrary state (probability distribution on $E$) $x\in S^{m-1}$
then passing to the state $V(x)$ (in the next "generation"), then to
the state $V(V(x))= V^2(x)$, and so on.

For a given $x^{(0)}\in S^{m-1}$ the trajectory (orbit)
$$ \{ x^{(n) }\}, \ \ n=0,1,2,...\ \ \mbox{of}\ \ x^{(0)}$$
under the action of QSO (3) is defined by
$$x^{(n+1)}= V(x^{(n)}),\ \ \mbox{where} \ \ n= 0,1,2,...$$

One of the main problem in mathematical biology consists in the study of
 the asymptotical behavior of the trajectories. The difficulty of the problem
 depends on given matrix \textbf{P}.

 For the history of (particularly)
 studied QSOs see [10], [11].

{\it The Volterra operators.}(see [3]-[5]) A Volterra QSO is defined by (3), (4)
 and the additional assumption
$$ P_{ij,k}=0, \ \ \mbox{if}\ \ k\not\in \{i,j\},\,\,
\forall i,j,k\in E. \eqno(5)$$

The biological treatment of condition (5) is clear: The offspring
repeats the genotype of one of its parents.

In paper [3] the general form of Volterra QSO
$$ V: x= (x_{1},...,x_{m}) \in S^{m-1}\,
\rightarrow\,V(x)= x'= (x'_1,...,x'_m)\,\in S^{m-1}
$$ is given
$$x'_k=x_k\left(1+\sum^m_{i=1}a_{ki}x_i\right),\eqno(6)$$
where
$$a_{ki}=2P_{ik,k}-1 \ \ \mbox{for}\, i\neq k\,\,\mbox{and}\
\,a_{ii}=0, i\in E.$$  Moreover
$$a_{ki}=-a_{ik}\ \ \mbox{and} \ \ |a_{ki}| \leq 1.$$

  In [3], [4] the theory of QSO (6) was developed
by using theory of the Lyapunov function and tournaments. But
non-Volterra QSOs (i.e. which do not satisfy the condition (5)) were
not in completely studied. Because there is no any general theory
which can be applied for investigation of non-Volterra operators.

In this paper we consider the following class of non-Volterra operators.

{\it $\ell$-Volterra QSO.} Fix $\ell \in E$ and assume that elements $P_{ij,k}$ of the matrix
$\textbf{P}$ satisfy
$$ P_{ij,k}= 0 \ \ \mbox{if} \ \ k \not\in
\{i,j\}\ \ \ \mbox{for}\ \ \mbox{any}\ \ k\in \{1,...,\ell\},\ \ i,j \in E;\eqno(7)$$
$$ P_{ij,k}> 0 \ \ \mbox{for at least one pair} \ \ (i,j),\ \
i\neq k,\ \ j\neq k \ \ \mbox{if} \ \  k \in \{\ell+1,...,m\}.\eqno(8)$$

\textbf{Definition 1.} For any fixed $\ell \in E$, the QSO defined by (3), (4), (7) and (8) is
called $\ell$-Volterra QSO.

Denote by $\mathcal{V}_{\ell}$ the set of all $\ell$-Volterra QSOs.

\textbf{Remarks.} 1. The condition (8) guarantees that $\mathcal{V} _{\ell_{1}} \bigcap
\mathcal{V}_{\ell_{2}} = \emptyset $ for any $\ell_{1}\neq \ell _{2}$.

2. Note that $\ell$-Volterra QSO is Volterra if and only if $\ell=
m$.

4. The class of $\ell$-Volterra QSO for a given $\ell$ does not coincide with a class of
non-Volterra QSOs mentioned in [10],[11].

\section{Lyapunov functions of $\ell$-Volterra QSO.}

Let $k\in \{1,...,\ell\}$ then $P_{kk,i}=0$ for $i\ne k$ and
$$\sum_{i=1}^mP_{kk,i}=P_{kk,k}+\sum_{i=\ell+1}^mP_{kk,i}=1.$$
Using $P_{ij,k}=P_{ji,k}$ we get for  $k=1,...,\ell$
 $$x'_k=x_k\left(P_{kk,k}x_k + 2\sum^m_{i=1\atop i\neq
k}P_{ik,k}x_i\right)=x_k\left(1+(P_{kk,k}-1)x_k +\sum^m_{i=1\atop i\neq
k}(2P_{ik,k}-1)x_i\right).$$

Denote $a_{ki}=2P_{ik,k}-1, k\ne i$ and  $a_{kk}=P_{kk,k}-1$ then we obtain

$$V: \left\{\begin{array}{llllll}
x_k' =x_k\left(1+\sum^m_{i=1}a_{ki}x_i\right),\ \ k=
1,...,\ell\\[3mm]
x_k' =x_k\left(1+\sum^m_{i=1}a_{ki}x_i\right)+\sum^m_{{i,j= 1\atop i\neq k}\atop j\neq
k}P_{ij,k}x_ix_j, \ \ k= \ell+1,...,m.
\end{array}\right.\eqno (9)$$

Note that
$$a_{kk}\in [-1,0]; \ |a_{ki}|\leq 1; \ a_{ki}+a_{ik}=2(P_{ik,i}+P_{ik,k})-2\leq 0,\ \
i,k\in E.\eqno(10)$$

 \vskip 0.2 truecm

{\bf Lemma 1.} {\it If $\exists k_0\in \{1,...,\ell\}$ and $\delta\leq 0$ such that $a_{k_0i}\leq
\delta$ for any $i\in E$ then}
$$P_\delta=\left\{p=(p_1,...,p_{\ell})\in S^{\ell-1}: \sum^\ell_{k=1}a_{ki}p_k\leq \delta, \
\mbox{for any} \ \ i\in E\right\}\ne \emptyset.$$

{\bf Proof.} It is easy to see that $e^{(k_0)}=(0,...,0,e^{(k_0)}_{k_0}=1,0,...,0)\in P_\delta.$
Thus for sufficiently small $\varepsilon$ we have
$$\left\{p\in S^{\ell-1}: \|p-e^{(k_0)}\|<\varepsilon\right\}\subset
P_\delta,$$ where $\|p-e^{(k_0)}\|=\max_i|p_i-e_i^{(k_0)}|.$ Indeed $\varepsilon$ can be chosen as
follows. Take $p=(p_1,...,p_\ell)$ with $p_{k_0}=1-\varepsilon, \ \sum_{i=1\atop i\ne k_0}^\ell
p_i=\varepsilon$ then $\|p-e^{(k_0)}\|\leq \varepsilon$ and
$$\sum_{k=1}^\ell a_{ki}p_k\leq \left\{\max_{k\in \{1,...,\ell\}\atop k\ne
k_0}a_{ki}\right\}\sum_{k =1\atop k\ne k_0}^\ell p_k+a_{k_0i}(1-\varepsilon)\leq$$
$$\left\{\max_{k\in \{1,...,\ell\}\atop k\ne k_0}\{a_{ki},
0\}-a_{k_0i}\right\}\varepsilon+a_{k_0i}\leq \delta$$ for any $i\in E$ if
$$\varepsilon\leq \min_{i\in E}{\delta-a_{k_0i}\over \max_{k\in \{1,...,\ell\}\atop k\ne k_0}\{a_{ki},
0\}-a_{k_0i}}.$$ This completes the proof.

\vskip 0.2 truecm

Let $\{x^{(n)}\}_{n=1}^\infty$ be the trajectory of the point $x^0\in S^{m-1}$ under operator (9).
Denote by $\omega(x^0)$ the set of limit points of the trajectory. Since $\{x^{(n)}\}\subset
S^{m-1}$ and $S^{m-1}$ is compact, it follows that $\omega(x^0)\ne \emptyset.$ Obviously, if
$\omega(x^0)$ consists of a single point, then the trajectory converges, and $\omega(x^0)$ is a
fixed point of (9). However, looking ahead, we remark that convergence of the trajectories is not
the typical case for the dynamical systems (9). Therefore, it is of particular interest to obtain
an upper bound for $\omega(x^0)$, i.e., to determine a sufficiently "small" set containing
$\omega(x^0)$.

Denote $${\rm int}S^{m-1}=\{x\in S^{m-1}: \prod_{i=1}^m x_i>0\}.$$
 \vskip 0.2 truecm

{\bf Definition 2.} A continuous function $\varphi:{\rm int}S^{m-1}\to R$ is called a Lyapunov
function for the dynamical system (9) if the limit $\lim_{n\to \infty}\varphi(x^{(n)})$ exists for
any initial point $x^0$.

Obviously, if $\lim_{n\to \infty}\varphi(x^{(n)})=c$, then $\omega(x^0)\subset \varphi^{-1}(c).$
Consequently, for an upper estimate of $\omega(x^0)$ we should construct Lyapunov functions that
are as large as possible. \vskip 0.2 truecm

{\bf Theorem 1.} {\it If $\exists k_0\in \{1,...,\ell\}$ such that $a_{k_0i}\leq 0$ for any $i\in
E$ and $p=(p_1,...,p_\ell)\in P_0$ then $\varphi_p(x)=x_1^{p_1}...x_\ell^{p_\ell}$ is a Lyapunov
function for (9).}

{\bf Proof.} Compute $\varphi_p(Vx)$:
$$\varphi_p(Vx)=\prod_{k=1}^\ell x_k^{p_k}\left(1+\sum_{i=1}^ma_{ki}x_i\right)^{p_k}=
\varphi_p(x)\prod_{k=1}^\ell \left(1+\sum_{i=1}^ma_{ki}x_i\right)^{p_k}.$$ Using Young's inequality
$$b_1^{p_1}...b_\ell^{p_\ell}\leq p_1b_1+...+p_\ell b_\ell,$$
where $b_i>0$ and $p_i\geq 0, \sum_{i=1}^\ell p_i=1$ we obtain

$$\varphi_p(Vx)\leq \varphi_p(x)\left(1+\sum_{i=1}^m\left(\sum_{k=1}^\ell a_{ki}p_k\right)x_i\right).$$

Since $p\in P_0$ we have $\sum_{k=1}^\ell a_{ki}p_k\leq 0$ for any $i\in E.$ Consequently
$\varphi(Vx)\leq \varphi(x)$. This completes the proof.

 \vskip 0.2 truecm

{\bf Theorem 2.} {\it If $a_{ki}<0$ for any $k=1,...,r, \ (r\leq \ell)$ and $i=r+1,...,m$ then
 $$\varphi(x)=x_1+...+x_r, \ x=(x_1,...,x_m)\in {\rm int}S^{m-1}$$
 is a Lyapunov function for (9). Moreover $\sum_{n=0}^\infty\varphi(x^{(n)})<+\infty$
 for any trajectory $\{x^{(n)}\}$.}

 {\bf Proof.} Compute $\varphi(Vx)$ according (9):
 $$\varphi(Vx)=\sum_{k=1}^rx'_k=\varphi(x)+\sum_{k=1}^r\left(\sum_{i=1}^ma_{ki}x_i\right)x_k.\eqno(11)$$
 Since $a_{kk}\in [-1,0], \ a_{ki}+a_{ik}\leq 0$ (see (10)), it follows that
 $$\sum_{k=1}^r\sum_{i=1}^ra_{ki}x_kx_i=\sum^r_{k=1}a_{kk}x_k^2+\sum_{1\leq k<i\leq
 r}(a_{ki}+a_{ik})x_kx_i\leq 0.$$

 Therefore, by (11) we have
$$\varphi(Vx)\leq\varphi(x)+\sum_{k=1}^r\left(\sum_{i=r+1}^ma_{ki}x_i\right)x_k.\eqno(12)$$
Let $$\alpha=\min_{k\in\{1,...,r\}\atop i\in\{r+1,...,m\}}\{-a_{ki}\},$$ since $\alpha>0,$ (12)
gives us
$$\varphi(Vx)\leq \varphi(x)-\alpha\sum_{k=1}^r\left(\sum_{i=r+1}^mx_i\right)x_k
=\varphi(x)[1-\alpha+\alpha\varphi(x)].\eqno(13)$$

For any $x^0\in {\rm int}S^{m-1}$ we have $\varphi(x^0)<1$. Since $\alpha\leq 1$, it follows
$1-\alpha+\alpha\varphi(x^0)<1$. Therefore, it follows from (13) that the inequality
$$\varphi(x^{(n+1)})\leq \varphi(x^{(n)})[1-\alpha+\alpha\varphi(x^{(n)})]\leq
\varphi(x^0)[1-\alpha+\alpha\varphi(x^0)]^n$$ holds along the trajectory $\{x^{(n)}\}$. Thus,
$\varphi(x^{(n)})\to 0$ and also $\sum_{n=0}^\infty\varphi(x^{(n)})<+\infty$. Theorem is proved.

It is known that if $a_n$ and $b_n$ are two bounded sequences of nonnegative numbers and if
$a_{n+1}\leq a_n+b_n$, $n=1,2,...$ then it follows from $b_n\to 0$ that $\{a_n\}$ is dense in
$[\underline{\lim}a_n, \overline{\lim}a_n]$. Moreover, if $\sum_{n=1}^\infty b_n<+\infty$, then
$\lim_{n\to\infty}a_n$ exists.

Below we use this fact to construct new Lyapunov functions.

\vskip 0.2 truecm

{\bf Theorem 3.} {\it If conditions of theorem 2 are satisfied then
$$\psi_p(x)=x_1^{p_1}...x_r^{p_r}, \ r\leq \ell, x=(x_1,...,x_m)\in{\rm int}S^{m-1}$$
is a Lyapunov function of (9) for any} $p=(p_1,...,p_r)\in S^{r-1}$.

{\bf Proof.} Using Young's inequality, we get
$$\psi_p(x')\leq \psi_p(x)\sum_{k=1}^r\left(1+\sum_{i=1}^ma_{ki}x_i\right)p_k=
\psi_p(x)\left(1+\sum_{k=1}^r\left(\sum_{i=1}^ma_{ki}x_i\right)p_k\right).\eqno(14)$$ By conditions
we have
$$\sum_{k=1}^r\sum_{i=1}^ra_{ki}p_kx_i\leq \sum_{k=1}^rx_i; \ \ \sum_{k=1}^r\sum_{i=r+1}^ma_{ki}p_kx_i\leq
0.$$ Hence by (14) we get
$$\psi_p(x')\leq \psi_p(x)(1+\sum_{k=1}^rx_k).$$
Consequently, along any trajectory $\{x^{(n)}\}$ we have
$$\psi_p(x^{(n+1)})\leq \psi_p(x^{(n)})(1+\varphi(x^{(n)})),\eqno(15)$$
where $\varphi(x^{(n)})=\sum_{k=1}^rx_k^{(n)}.$ According to Theorem 2, the series
$\sum^\infty_{n=0}\varphi(x^{(n)})$ converges, and so it follows from (15) that
$\lim_{n\to\infty}\psi_p(x^{(n)})$ exists along any trajectory.

{\bf Remark.}  When the functions $\varphi_p$ and $\psi_p$ are extended from ${\rm int}S^{m-1}$ to
$S^{m-1}$ the expression $0^0$ can arise, and we set it equal to 1.

Now we shall describe Lyapunov functions of other forms.

\vskip 0.2 truecm

{\bf Theorem 4.} {\it If there exists $p\in \{1,...,\ell\}$ and $q\in E$ such that
$a_{pi}-a_{qi}\leq 0$ for any $i\in E$ then
$$f_{pq}(x)={x_p\over x_q}, \  x=(x_1,...,x_m)\in{\rm int}S^{m-1}$$
is a Lyapunov functions of (9). Moreover $f_{pq}(x)$ is monotonically decreasing along the
trajectory $\{x^{(n)}\}$, where $x^0\in {\rm int}S^{m-1}$ and $x^0\ne V(x^0)$.}

{\bf Proof.} We have
$$f_{pq}(x')={x'_p\over x'_q}=f_{pq}(x)\cdot {1+\sum_{i=1}^m a_{pi}x_i\over 1+\sum_{i=1}^m a_{qi}x_i+{\bf 1}_{\{q>\ell\}}
x^{-1}_q\sum_{i,j=1\atop i\ne q, j\ne q}^mP_{ij,q}x_ix_j},\eqno(16)$$ where ${\bf
1}_{\{q>\ell\}}=0$ (resp. =1) if $q\leq \ell$ (resp. $q> \ell$). Clearly,
$${\bf 1}_{\{q>\ell\}} x^{-1}_q\sum_{i,j=1\atop i\ne q, j\ne q}^mP_{ij,q}x_ix_j\geq 0, \
\mbox{for any} \ x\in {\rm int} S^{m-1}.$$ Consequently, from (16) by condition of theorem we get
$$f_{pq}(x')\leq \alpha f_{pq}(x),\eqno(17)$$
where $$\alpha=\max_{x\in {\rm int}S^{m-1}}{1+\sum_{i=1}^m a_{pi}x_i\over 1+\sum_{i=1}^m
a_{qi}x_i}\leq 1.$$ This implies $f_{pq}(x^{(n+1)})< f_{pq}(x^{(n)}), n\geq 0.$ Thus sequence
$f_{pq}(x^{(n)})$ is a monotonically decreasing. Since it is bounded we conclude that $f_{pg}$ is a
Lyapunov function.

{\bf Remark.} The Lyapunov functions mentioned in Theorems 1, 2 and 4 are monotonically decreasing
along any trajectory. Note that under conditions of Theorem 4 we can also construct the function
$f^+_{qp}(x)={x_q\over x_p}$ which is monotonically increasing along any trajectory $\{x^{(n)}\}$
with $x^{0}\in {\rm int}S^{m-1}$. But the limit $\lim_{n\to \infty}f^+_{qp}(x^{(n)})$ can be equal
to $+\infty$.

\section{Upper estimations of $\omega(x^0)$}

In this section we shall apply the Lyapunov functions described in the previous section to obtain
an upper bound of $\omega(x^0)$.

Denote by ${\rm Fix}(V)$ the set of all fixed points of the operator (9) i.e.
$${\rm Fix}(V)=\{x\in S^{m-1}: V(x)=x\}.$$

 \vskip 0.2 truecm
{\bf Theorem 5.} {\it If there exists $k_0\in \{1,...,\ell\}$ and $\delta>0$ such that
$a_{k_0i}\leq -\delta$ for any $i\in E$ then for} $x^0\notin {\rm Fix}(V)$,
$$\omega(x^0)\subset \{x\in S^{m-1}: \prod_{i=1}^\ell x_i=0\}.$$

{\bf Proof.} Consider Lyapunov function $\varphi_p(x)=\prod_{i=1}^\ell x_i^{p_i}$ for $p\in
P_\delta$. By proof of Theorem 1 we have
$$\varphi_p(Vx)\leq (1-\delta)\varphi_p(x), \ \delta>0.$$
Iterating this inequality we obtain $\varphi_p(x^{(n)})\leq (1-\delta)^n\varphi_p(x^0).$ Hence
$$\lim_{n\to\infty}\varphi_p(x^{(n)})=\lim_{n\to\infty}\prod_{i=1}^\ell(x_i^{(n)})^{p_i}=0.$$ This
completes the proof.

\vskip 0.2 truecm

As a corollary of Theorem 2  we have

{\bf Theorem 6.} {\it Suppose conditions of Theorem 2 are satisfied. If $i\in \{1,...,r\}$ then
$x_i^{(n)}\to 0$, at the rate of a geometric progression as $n\to \infty$.}

This Theorem gives the estimation $\omega(x^0)\subset S^{m-r-1}$, where $S^{m-r-1}$ is the face of
$S^{m-1}$ spanned by the vertices $e^{(r+1)}=(0,...,0,e_{r+1}=1,0,...,0)$,...,
$e^{(m)}=(0,...,0,1)$.

If in Theorem 4 we consider more stronger condition i.e. $a_{pi}-a_{qi}<0$ instead of
$a_{pi}-a_{qi}\leq 0$, for any $i\in E$. Then we get (17) with $\alpha<1$. In this case it follows
that $f_{pq}(x^{(n)})\to 0$. Using the fact that $0<x_q^{(n)}<1$, we get $x_p^{(n)}\to 0$. This
enables us to get a more precise estimate for $\omega(x^0):$ it is a subset of the simplex
$S^{m-1}$ with $x_p=0$ where $p\in \{1,...,\ell\}$ such that there exists $q=q(p)\in E$ which
satisfies conditions of Theorem 4.

By these results and results of [11] we make following remarks.\\

 {\bf Remarks.} 1. For Volterra operators the estimate
$\omega(x^0)\subset S^{m-r-1}$ can be improved to the estimation $\omega(x^0)\subset \partial
S^{m-r-1}=\{x\in S^{m-r-1}: \prod_{i=1}^{m-r-1}x_i=0\}$ (see [3]). In general, if $\ell<m$ then
such an improvement is impossible.

2. If $\ell\leq m-2$ then $\ell$-Volterra operators can have cyclic trajectories this is quite
different behavior from the behavior of Volterra operators, since Volterra operators have no cyclic
trajectories.

3. One of the main goal by introducing the notion of
 $\ell$-Volterra operators was to give an example of QSO which has
 more rich dynamics than Volterra QSO. It is well known [3], [5] that
 for Volterra operators (see (6)) if $a_{ij}\neq 0 \ \ (i\neq j)$ then
 for  any non-fixed initial point $\lambda^{0}$ the set
 $\omega(\lambda^{0})$ of all limit points of the trajectory
 $\{\lambda^{(n)}\}$ is subset of the boundary of simplex. But for $\ell$-Volterra
 operators, in general, the limit set can be  subset of the inside of simplex.

4. It is known [3] that Volterra operators are homomorphisms. Consequently, for any initial point
$x^0\in S^{m-1}$ the "negative" trajectory $\{V^{-n}(x^0)\}, n=0,1,2,...$ exists. Moreover the
negative trajectories always converge. But such kind of result is not true for $\ell$-Volterra
operators.

\section{The fixed points of the operator (9)}

It is easy to see that a vertex $e^{(i)}=(0,...,0,1_i,0,...,0)$  of $S^{m-1}$ is a fixed point of
$V$ iff $P_{ii,i}=1$. We consider the question of the existence of other fixed points.

For $j\in E$ denote
$$X_{j}=\left\{x\in S^{m-1}: x'_k=V(x)_k=x_k, \ \ k=1,...,j\right\}.$$

Note that $X_{m-1}=X_m={\rm Fix}(V)$ and $X_j\subset X_{j-1}$ for any $j=1,...,m-1$.

It is easy to see that $x\in X_\ell$ if and only if
$$x_k\sum^m_{i=1}a_{ki}x_i=0, \ k=1,...,\ell.\eqno(18)$$

Set ${\rm supp}_\ell x=\{i\in \{1,...,\ell\}: x_i\ne 0\}$ then from (18) we get

$${\rm supp}_\ell x\cap {\rm supp}_\ell Ax=\emptyset,\eqno(19)$$
where $A=(a_{ij})_{i,j=1}^m$ is $m\times m$ matrix with $a_{ij}$ defined in section 3.

\vskip 0.2 truecm

{\bf Lemma 2.} {\it If $x,y\in X_\ell$ and ${\rm supp}_\ell x={\rm supp}_\ell y$ then $\lambda
x+(1-\lambda)y\in X_\ell$ for any $\lambda \in [0,1]$.}

{\bf Proof.} By (19) and ${\rm supp}_\ell x={\rm supp}_\ell y$ we have
$${\rm supp}_\ell x\cap({\rm supp}_\ell Ax\cup {\rm supp}_\ell Ay)=\emptyset.$$
Since ${\rm supp}_\ell (\lambda u+(1-\lambda)v)\subset {\rm supp}_\ell u\cup {\rm supp}_\ell v$ we
have
$${\rm supp}_\ell (\lambda x+(1-\lambda)y)\cap {\rm supp}_\ell(\lambda Ax+(1-\lambda)Ay)\subset
{\rm supp}_\ell x\cap ({\rm supp}_\ell Ax\cup {\rm supp}_\ell Ay)=\emptyset.$$ Hence $\lambda
x+(1-\lambda)y\in X_\ell$.

\vskip 0.2 truecm For $p,q,r\in E$ put
$\Delta=(a_{pr}-a_{pp})(a_{qr}-a_{qq})-(a_{pr}-a_{pq})(a_{qr}-a_{qp});$

$$\Delta_1=a_{qr}a_{pq}-a_{pr}a_{qq}; \ \Delta_2=a_{pr}a_{qp}-a_{qr}a_{pp}.$$
\vskip 0.3 truecm
 {\bf Theorem 7.} {\it If}

 (a) $p,q,r\in E$ {\it such that $\{p,q,r\}\cap \{\ell+1,...,m\}$ contains at most one
element, say $r$;}

(b) $P_{ij,k}=0$ {\it for any} $i,j\in \{p,q,r\}, k\in \{\ell+1,...,m\}\setminus \{p,q,r\}$;

(c) $\Delta\ne 0$, ${\rm sign}(\Delta)= {\rm sign}(\Delta_1)={\rm sign}(\Delta_2)={\rm
sign}(\Delta-\Delta_1-\Delta_2).$

{\it Then the interior of the two-dimensional face $S^2_{pqr}$ spanned by the vertices $e^{(p)},
e^{(q)}$ and $e^{(r)}$ of $S^{m-1}$ contains exactly one fixed point of $V$.}

{\bf Proof.} Under conditions (a) and (b) the restriction of V (see (9)) to the face $S^2_{pqr}$
has the form

$$\begin{array}{lll}
    x'_p=x_p\left(1+a_{pp}x_p+a_{pq}x_q+a_{pr}x_r\right), \\
   x'_q=x_q\left(1+a_{qp}x_p+a_{qq}x_q+a_{qr}x_r\right), \\
   x'_r=x_r\left(1+a_{rp}x_p+a_{rq}x_q+a_{rr}x_r\right)+P_{pp,r}x_p^2+2P_{pq,r}x_px_q+P_{qq,r}x_q^2.\\
\end{array}\eqno(20)
$$
Using $x_p+x_q+x_r=1$, $x_px_qx_r>0$ and (20) the equation $V(x)=x$ on $S^2_{pqr}$ can be written
as
$$\begin{array}{ll}
  \left(a_{pr}-a_{pp}\right)x_p+\left(a_{pr}-a_{pq}\right)x_q=a_{pr},\\
  \left(a_{qr}-a_{qp}\right)x_p+\left(a_{qr}-a_{qq}\right)x_q=a_{qr}. \\
\end{array}
$$
Elementary computations show that $x^*=(x^*_1,...,x^*_m)$, where
$$x^*_p={\Delta_1\over \Delta},\ \  x^*_q={\Delta_2\over \Delta},\ \ x^*_r={\Delta-\Delta_1-\Delta_2\over
\Delta}, $$ and all the rest of the coordinates are zero, is a fixed point of $V$. It follows from
condition (c) that $x^*_p, x^*_q, x^*_r>0$. Therefore, $x^*$ is a fixed point satisfying the
condition of the theorem. Uniqueness is verified by a simple computation. Theorem is proved. \vskip
0.2 truecm

{\bf Acknowledgments.} This work was done within the scheme of Junior Associate at the ICTP,
Trieste, Italy and UAR thanks ICTP for providing financial support and all facilities (in June -
August 2008). He also thanks Prof. S.Albeverio for an invitation to Bonn University (Germany) in
August 2008 and for many helpful discussions.

\vskip 0.2 truecm {\bf References}

1. Devaney R. L., An introduction to chaotic dynamical system, {\it Westview Press,} (2003).

2. Elaydi S. N., Discrete chaos, \textit{Chapman  Hall/CRC,} (2000).

3. Ganikhodzhaev R.N., Quadratic stochastic operators, Lyapunov functions and tournaments, {\it
Russian Acad. Sci. Sbornik Math.}, {\bf 76} (1993), 489-506.

4. Ganikhodzhaev R.N., A chart of fixed points and Lyapunov functions for a class of discrete
dynamical systems. {\it Math. Notes}, {\bf 56} (1994), 1125-1131.

5. Ganikhodzhaev R. N., Eshmamatova D.B., Quadratic automorphisms of simplex and asymptotical
behavior of their trajectories. {\it Vladikavkaz Math. Jour.}, {\bf 8} (2006), 12-28.

6. Hofbaver J., Sigmund K., The theory of evolution and dynamical systems, \textit{Cambridge Univ.
Press}, (1988).

7. Kesten H, Quadratic transformations: A model for population growth I, II. {\it Adv. Appl.
Prob.}, {\bf 2} (1970), 1-82 and 179-228.

8. Lyubich Yu. I., Mathematical structures in population genetics, {\sl Biomathematics},
Springer-Verlag, {\bf 22} (1992).

9. Robinson R. C, An introduction to Dynamical systems: Continues and Discrete, \textit{Pearson
Education}, (2004).

10. Rozikov U. A., Jamilov U. U., On F-Quadratic stochastic operators,  {\it Math. Notes}, {\bf 83}
(2008), 606-612.

11. Rozikov U. A., Zada A., On Dynamics of $\ell$- Volterra Quadratic Stochastic Operators,
arXiv:0712.3979 [math.DS].

\end{document}